\documentclass[paper=a5,11pt,headings=small]{scrartcl}

\usepackage[english]{babel}
\usepackage[T1]{fontenc}
\usepackage[utf8]{inputenc}
\usepackage{mathptmx}
\usepackage[left=17mm,right=20mm,top=20mm,bottom=20mm,            footskip=2\baselineskip]{geometry}

\usepackage[outline]{contour}
\usepackage{overpic}
\usepackage{amsmath,amsfonts,amssymb,amsthm}
\usepackage{mathtools}
\usepackage{algorithm}
\usepackage{algpseudocode}

\usepackage{microtype}
\usepackage[sort&compress]{natbib}
\setcitestyle{square,numbers}

\usepackage{hyperref}
\hypersetup{
  pdfauthor={Gábor Hegedüs, Zijia Li, Josef Schicho, Hans-Peter Schröcker},  pdftitle={From the Fundamental Theorem of Algebra to Kempe's Universality Theorem}}

\usepackage{scrpage2}
\addtokomafont{title}{\flushleft}
\setkomafont{pagefoot}{\normalfont}
\ofoot{\thepage}
\title{From the Fundamental Theorem of Algebra to Kempe's Universality Theorem}
\author{  Gábor Hegedüs\thanks{Applied Mathematical Institute, Antal Bejczy Center for Intelligent Robotics, Obuda University, 1032 Budapest, Hungary} \and
  Zijia Li\thanks{Johan Radon Institute for Computational and Applied Mathematics, Austrian Academy of Sciences, 4040 Linz, Austria} \and
  Josef Schicho\thanks{Research Institute for Symbolic Computation, Johannes Kepler University Linz, Schloss Hagenberg, 4232 Hagenberg, Austria} \and
  Hans-Peter Schröcker\thanks{Unit Geometry and CAD, University of Innsbruck, 6020 Innsbruck, Austria}}
            
\setkomafont{sectioning}{\normalfont\bfseries}

\newcommand{\D}{\mathbb{D}}
\renewcommand{\H}{\mathbb{H}}
\renewcommand{\DH}{\D\H}
\newcommand{\R}{\mathbb{R}}
\newcommand{\SSS}{\mathbb{S}}
\newcommand{\eps}{\varepsilon}
\newcommand{\qi}{\mathbf{i}}
\newcommand{\qj}{\mathbf{j}}
\newcommand{\qk}{\mathbf{k}}
\newcommand{\cj}[1]{\overline{#1}}
\newcommand{\SE}[1][3]{\mathrm{SE}(#1)}
\newcommand{\SQ}{\mathcal{S}}
\newcommand{\EG}{E}

\newtheorem{theorem}{Theorem}
\newtheorem{proposition}{Proposition}
\newtheorem{corollary}{Corollary}
\theoremstyle{remark}
\newtheorem{example}{Example}

\begin{document}

\maketitle

This article provides a gentle introduction for a general mathematical
audience to the factorization theory of motion polynomials and its
application in mechanism science. This theory connects in a rather
unexpected way a seemingly abstract mathematical topic, the non-unique
factorization of certain polynomials over the ring of dual
quaternions, with engineering applications. Four years after its
introduction \cite{hegedues12:_factorization,  hegedus13:_factorization2}, it is already clear how beneficial it
has been to both fields \cite{gallet15,  hegedus15:_four_pose_synthesis,  li15:_rational_motions,  li15:_motion_polynomials,  li15:_darboux_linkages,  li14:_sharp,  li14:_straight,  li13:_anglesymmetric,li14:_ck}. In
Section~\ref{sec:motion-polynomials} we introduce the notion of motion
polynomials and discuss their decomposition into products of linear
motion polynomials. It can be used to synthesize linkages following a
prescribed motion and is related to a variant of Kempe's Universality
Theorem. We explain the relation to mechanism science in more detail
in Section~\ref{sec:factorizations-linkages}. In
Sections~\ref{sec:synthesis} and \ref{sec:exceptional-factorizations}
we present examples from linkage synthesis and discuss exceptional
factorizations.

\section{Motion polynomials and their factorizations}
\label{sec:motion-polynomials}

We denote by $\H$ the skew field of quaternions, generated by
$\qi,\qj,\qk$ over $\R$, with the well-known relations
\begin{equation*}
  \qi^2 = \qj^2 = \qk^2 = \qi\qj\qk= -1.
\end{equation*}
The conjugate of a quaternion $q=q_1+q_2\qi+q_3\qj+q_4\qk$ is defined
by $\cj{q}=q_1-q_2i-q_3j-q_4k$, the norm of $q$ is
$N(q) \coloneqq q\cj{q} = q_1^2 + q_2^2 + q_3^2 + q_4^2$. The scalar
extension $\DH \coloneqq \D\otimes_\R \H$ by \emph{dual numbers}
$\D \coloneqq \R[\eps]/\langle\eps^2\rangle$ is just a skew ring, the ring of
\emph{dual quaternions}.  The conjugate of a dual quaternion
$h=p+\eps q$ is defined by $\cj{h}=\cj{p}+\eps\cj{q}$. We also use
$N(h)$ to denote the norm of $h$ with
$N(h) \coloneqq h\cj{h} = N(p) + \eps(p\cj{q} + q\cj{p}) \in \D$.  The
dual quaternion $h$ is invertible if and only if $p \neq 0$.

Denote by $\SSS$ the multiplicative subgroup of dual quaternions with
nonzero real norm. Its elements may be written as $h = p + \eps q$
where the quaternions $p$ and $q$ satisfy $p \neq 0$,
$p\cj{q} + q\cj{p} = 0$. The latter equality is called the \emph{Study
  condition.} The group $\SSS$ acts on
$\R^3 = \langle \qi, \qj, \qk \rangle$ according to
\begin{equation}
  \label{eq:3}
  x \mapsto \frac{px\cj{p}+p\cj{q} - q\cj{p}}{N(p)}.
  \end{equation}
Any such map is a \emph{Euclidean displacement}. The rotational part
is the well-known action $x \mapsto px\cj{p}/N(p)$ of the unit
quaternion $p N(p)^{-1/2}$ on $\R^3$, the translational component is
$(p\cj{q} - q\cj{p})/N(p)$. Equation~\eqref{eq:3} defines a
\emph{homomorphism} from $\SSS$ to the group $\SE$ of Euclidean
displacements. It is surjective and its kernel is the real
multiplicative group~$\R^\ast$.

This algebraic construction allows a geometric interpretation.
Identify $\DH / \R^\ast$ with real projective space $P^7$ of dimension
seven, denote by $\SQ \subset P^7\colon p\cj{q} + q\cj{p} = 0$ the
\emph{Study quadric} and by $\EG$ the three-space with equation
$p = 0$. Then \emph{Study's kinematic map} is the bijection
$\SE \cong \SSS / R^\ast \to \SQ \setminus \EG \subset P^7$ whose
inverse maps $h = p + \eps q$ to the Euclidean displacement given by
\eqref{eq:3}.

The concept of motion polynomials arises by making the above group
homomorphism parametric. More precisely, let $\DH[t]$ be the skew ring
of univariate polynomials over $\DH$ in one variable $t$ that commutes
with all coefficients. For $C \in \DH[t]$, the conjugate polynomial
$\cj{C}$ is obtained by conjugating all coefficients.  If $C=P+\eps Q$
with $P,Q\in\H[t]$, then we call $P$ the primal part and $Q$ the dual
part. We say that $C\in\DH[t]$ is a \emph{motion polynomial} if
$N(C) \coloneqq C\cj{C} \in \R[t] \setminus \{0\}$ (a priori, the norm
polynomial $N(C)$ is in $\D[t]$) and if its leading coefficient is
invertible. For any $t_0\in\R$ which is not a zero of $N(C)$, we can
say that $C(t_0)$ is an element of $\SSS$, hence it acts on
$\R^3$. Varying $t_0$, we get a motion, i.e. a parametrized curve of
Euclidean displacements.

By virtue of \eqref{eq:3}, the trajectories of points during that
motion are rational curves whence the motion itself is called
\emph{rational.} It is well-known that all motions with only rational
trajectories have a polynomial parametrization with values in the
Study quadric \cite{juettler93:_rationale_bewegungsvorgaenge}. The
motion polynomials do not constitute a multiplicative group.  However,
the quotient by $\R[t]-\{0\}$ is a group because the inverse of the
class of a motion polynomial is precisely its conjugate. While
kinematic or geometric properties of motion polynomials do not change
if we multiply them with non-zero real polynomials, algebraic
properties may be different. This observation will be important in
Section~\ref{sec:exceptional-factorizations}.

Summarizing, we can state
\begin{proposition}
  Motion polynomials parametrize rational motions and the group of
  motion polynomials modulo $\R[t] \setminus \{0\}$ is isomorphic to
  the group of rational motions. Via Study's kinematic mapping, motion
  polynomials (rational motions) correspond to rational curves on the
  Study quadric $\SQ$ with at most finitely many points in the three
  space~$\EG$.
\end{proposition}

Since we want to use some version of the fundamental theorem of
algebra later on, it is convenient to restrict attention to monic
motion polynomials. As far as applications in kinematics are
concerned, this is no loss of generality and can always be
accomplished by a suitable coordinate change.  The fundamental theorem
speaks about factorization into linear polynomials, so it is time to
clarify the kinematic nature of linear motion polynomials.

\begin{proposition}
  \label{prop:3}
  Every monic linear motion polynomial parametrizes either a rotation
  about a fixed axis in $\R^3$ or a translation in a fixed direction
  in~$\R^3$.
\end{proposition}

The converse is not true: there are monic motion polynomials
parametrizing a rotation around a fixed axis that are not linear.

\begin{example}
  The dual quaternion polynomial $C_1=t-\qi$ has norm $t^2+1$, hence
  it is a motion polynomial.  It parametrizes a rotation around the
  first coordinate axis. Indeed, by \eqref{eq:3} we have
  \begin{multline*}
    x_1\qi + x_2\qj + x_3\qk \mapsto
    \frac{(t-\qi)(x_1\qi + x_2\qj + x_3\qk)(t+\qi)}{t^2+1} \\
    = x_1\qi + \Bigl(\frac{t^2-1}{t^2+1}x_2 + \frac{2t}{t^2+1}x_3\Bigr)\qj
             + \Bigl(\frac{t^2-1}{t^2+1}x_3 - \frac{2t}{t^2+1}x_2\Bigr)\qk.
  \end{multline*}
  Setting $\varphi \coloneqq 2 \arctan t$, i.e.,
  $\cos\varphi = (1-t^2)/(t^2+1)$, $\sin\varphi = 2t/(t^2+1)$, the
  assertion becomes apparent.

  The dual quaternion polynomial $C_2=(t-\qi)^2=t^2-1-2\qi t$ has norm
  $(t^2+1)^2$, hence it is a motion polynomial. It also
  parametrizes a rotation around the first coordinate axis. Indeed, it
  can be obtained by reparametrizing the first parametrization using
  the reparametrization function $t\mapsto (t^2-1)/(2t)$ and clearing
  denominators.
\end{example}

\begin{example}
  The dual quaternion polynomial $C_3=t-\eps\qi$ has norm $t^2$, hence
  it is a motion polynomial. The action on $\R^3$ is
  \begin{equation*}
    x_1\qi + x_2\qj + x_3\qk \mapsto
    x_1\qi + x_2\qj + x_3\qk + \frac{2\eps\qi}{t}.
  \end{equation*}
  This is a translation in the direction parallel to the first
  coordinate axis.
\end{example}

The following theorem is fundamental in the factorization theory
of quaternion polynomials \cite{gordon65}:
\begin{proposition}
  \label{prop:4}
  Every monic polynomial $C \in \H[t]$ admits a factorization
  $C = (t-q_1) \cdots (t-q_n)$ with quaternions $q_1,\ldots,q_n$.
\end{proposition}

For the proof, one uses a fundamental theorem for monic non-negative
real polynomials: any such polynomial (it must be of even degree) has
a unique factorization into monic non-negative real quadratic
polynomials. This is an easy consequence of the fundamental theorem
for complex polynomials.  The norm polynomial $N(C)$ is a non-negative
real polynomial and factors into the product of non-negative real
quadratic polynomials. The Euclidean algorithm for polynomial division
also works in the case of polynomials in $\H[t]$. If we take one of
the non-negative quadratic factors of the norm, say $Q \in \R[t]$,
then the polynomial remainder of $C$ by $Q$ is a linear quaternion
polynomial or zero. In the first case, this linear polynomial can be
factored out. In the second case, $Q$ can be factored out.

A specialty in this context is that a generic polynomial $C \in \H[t]$
has precisely $n!$ different factorizations. Since the ring
of quaternion polynomials is not commutative, permuting the factors no
longer yields a factorization $C$. In other words the $n!$ factorizations
differ essentially, whereas in the complex case all factorizations can
be obtained as permutations of one. In this case we rather would think of
just one factorization whose factors can be permuted.

Another word of caution concerns the polynomial division mentioned in
the above sketch of the proof. Since $\H[t]$ is not commutative, we
have left and right polynomial division, left and right quotients, and
left and right remainders. It does not matter which direction we
choose, but we have to choose one and stick to it.  In the following
we deal only with right polynomial division, right remainders, and
consequently left quotients:
($C = \text{Quotient}\cdot D + \text{Remainder}$).

Let us pass from polynomials over $\H[t]$ to motion polynomials in
$\DH[t]$ and try to factor the motion polynomial $C$. Like in the
proof above, the norm is non-negative and, by the important defining
property of motion polynomials, \emph{real.} The quotient by one of
the quadratic factors has degree at most one. But now, polynomial
division might be problematic because the linear remainder polynomial
is in general not monic and the leading coefficient might not be
invertible. It is even possible that the remainder is constant, as in
the example $C=t^2+1+\eps\qi$ where we have $C\cj{C} = (t^2+1)^2$ and
the remainder is~$\eps\qi$.

However, in the generic case, everything is fine: the remainders will
be linear with invertible leading coefficients and can be divided
out. We can then construct step by step a factorization into linear
factors by Algorithm~\ref{alg:factor} below. This genericity condition
can be simplified as follows: For Algorithm \ref{alg:factor} to
work, it is sufficient that the primal part of $C$ has no nontrivial
real factors. The different factorizations come from the arbitrary
choice of a quadratic factor $M$ in Line~\ref{line:choose} of the
algorithm.

\begin{algorithm}
  \caption{(factorization of generic motion polynomials)}
  \label{alg:factor}
  \begin{algorithmic}[1]
    \Require Motion polynomial $C$, monic, primal part has no real
    factor.    \Ensure List $L = [h_1,\ldots,h_n]$ of generic motion polynomials\newline
    such that $C = (t-h_1)\cdots(t-h_n)$.    \State \label{line:qfactors} $F \leftarrow$ list of quadratic,
    irreducible factors of $C\cj{C}$    \State $D \leftarrow C$ \State initialize $L = [\;]$ \Comment
    empty list \While{$F$ is not empty} \State \label{line:choose}
    choose $M \in F$    \State write $D = QM + R$ with $\deg R \le 1$ \Comment polynomial
    right division \State \label{line:zero} compute unique zero $h$ of
    $R$    \State prepend $h$ to $L$    \State $D \leftarrow D'$ where $D = D'(t-h)$ \Comment polynomial
    right division \State delete $M$ from $F$    \EndWhile
  \end{algorithmic}
\end{algorithm}

In fact, we have

\begin{theorem}[\cite{hegedus13:_factorization2}]
  \label{th:1}
  Algorithm~\ref{alg:factor} can be used to factor a motion polynomial
  $C = P + \eps Q$, provided the primal part $P$ has no real
  factors. In this case, $C = (t-h_1) \cdots (t-h_n)$ and for
  $i \in \{1,\ldots,n\}$ each polynomial $t - h_i$ describes a
  rotation about a fixed axes. Moreover, all possible factorizations
  (in general $n!$) can be obtained in that way.
\end{theorem}

As a consequence of Theorem~\ref{th:1} we have
\begin{corollary}
  \label{cor:1}
  In general, a rational motion of degree $n$ can be decomposed in at
  most $n!$ different ways into the product of rotations $t-h_1$,
  \ldots, $t-h_n$.
\end{corollary}

The phrase ``in general'' in Corollary~\ref{cor:1} refers to the absence of
real factors in the primal part of $C$. If this requirement is not
fulfilled, the statement is not true. In fact, all kinds of special
behavior can be observed. It is interesting that these ``exceptional''
examples surprisingly often arise in natural applications of
Theorem~\ref{th:1} to kinematics and mechanism science. We will return to
this point a little later. But at first, we have to explain the
relation between Theorem~\ref{th:1} and mechanism science.

\section{Factorizations and linkages}
\label{sec:factorizations-linkages}

\begin{figure}
  \centering
  \includegraphics{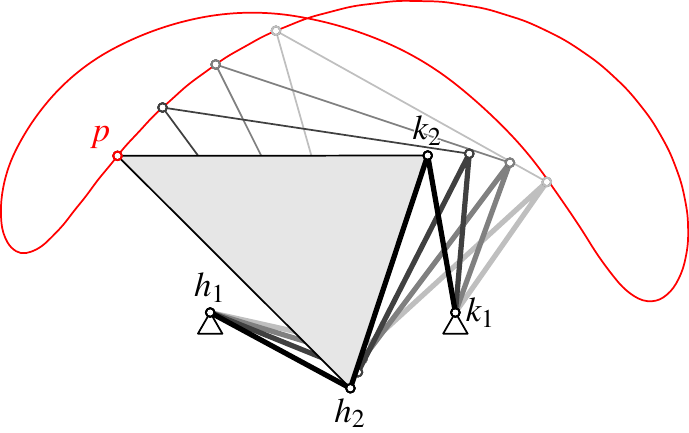}
  \caption{Anti-parallelogram linkage}
  \label{fig:antiparallelogram}
\end{figure}

Consider a generic motion polynomial $C$ of degree two that
parametrizes a planar motion (all trajectories are in the parallel
planes). By Theorem~\ref{th:1}, it admits two factorizations
\begin{equation*}
  C = (t - h_1)(t - h_2) = (t - k_1)(t - k_2)
\end{equation*}
with suitable dual quaternions $h_1,h_2,k_1,k_2$. This means that the
motion parametrized by $C$ can be generated in two ways as composition
of two rotations with respective centers $h_1$, $h_2$ or $k_1$, $k_2$.
Moreover, the two revolute joints at $h_2$ and $k_2$ can be rigidly
connected without disturbing the motion $C$. This is illustrated in
Figure~\ref{fig:antiparallelogram} which depicts a planar four-bar
linkage whose joints are labeled by the corresponding dual
quaternions. The joints $h_1$ and $k_1$ are fixed, $h_2$ and $k_2$ can
rotate about $h_1$ and $k_1$, respectively, but retain their
distance. Planar four-bar linkages constitute the most important class
of linkages for engineering applications. Our case is special because
the joints $h_1$, $h_2$, $k_2$, and $k_1$ form an anti-parallelogram.

If the motion polynomial $C$ is not planar, a similar construction
yields a spatial linkage, consisting of a closed loop of four skew
revolute axes (Figure~\ref{fig:bennett}). Any two neighboring axes are
rigidly connected, that is, they maintain their distance and angle
throughout the motion. In contrast to the planar case, the
one-dimensional mobility of such a structure is not obvious. A closed
loop of four revolute joints is generically rigid. But the loops
resulting from factorizations of quadratic motion polynomials move
because of their algebraic construction.

Spatial four-bar linkages that move with one degree of freedom have
been known for a long time \cite{bennett03,bennett14,krames37}.  There
exists only one family of this linkage type and its members are
referred to as ``Bennett linkages''.  So far, Bennett linkage have
minor importance in applications but we shall see their theoretical
significance later in this text. Moreover, they exhibit a fascinating
geometry. A stunning example is the Wunderlich's explanation of the
Bäcklund transform of discrete asymptotic nets of constant Gaussian
curvature by means of Bennett linkages
\cite{wunderlich51:_differenzengeometrie}.

By factorizing cubic motion polynomials, we can also construct spatial
six-bar linkages with a one-dimensional mobility as in
Figure~\ref{fig:6R}. Their complete classification is a long-standing
open problem in theoretical mechanism science. In spite of some recent
progress \cite{hegedus13:_bonds2,hegedus15:_maximal_genus}, a
classification is currently still out of reach. At any rate, our approach
yields new examples of spatial six-bar linkages
\cite{li13:_anglesymmetric,li14:_ck,li14:_sharp} and the first viable
synthesis procedure. We describe this in more details in the next
section.

\section{Linkage synthesis}
\label{sec:synthesis}

An important application of motion polynomial factorization is
\emph{linkage synthesis,} the construction of a mechanical linkage to
accommodate a certain task. Our approach is well suited to exact
synthesis with prescribed poses, that is, the computation of linkages
such that one link visits a finite number of prescribed poses. 
The word ``pose'' refers here to the position and orientation of a rigid
body, that is, an element of $\SE$.

\begin{figure}
  \centering
  \small
  \begin{overpic}{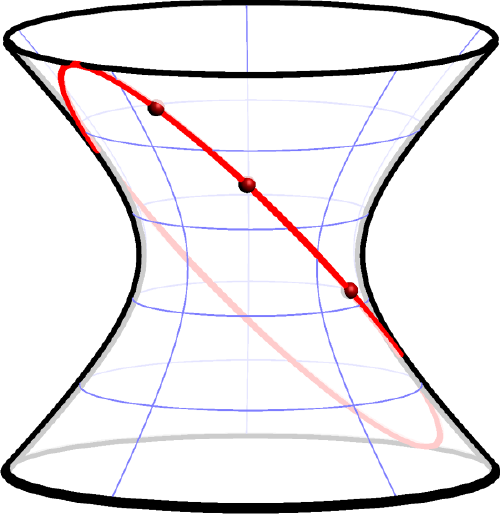}
    \put(26,71){\contour{white}{$p_0$}}
    \put(40,57){\contour{white}{$p_1$}}
    \put(61,36){\contour{white}{$p_2$}}
  \end{overpic}
  \hfill
  \begin{overpic}[page=2,trim=20 60 20 58,clip]{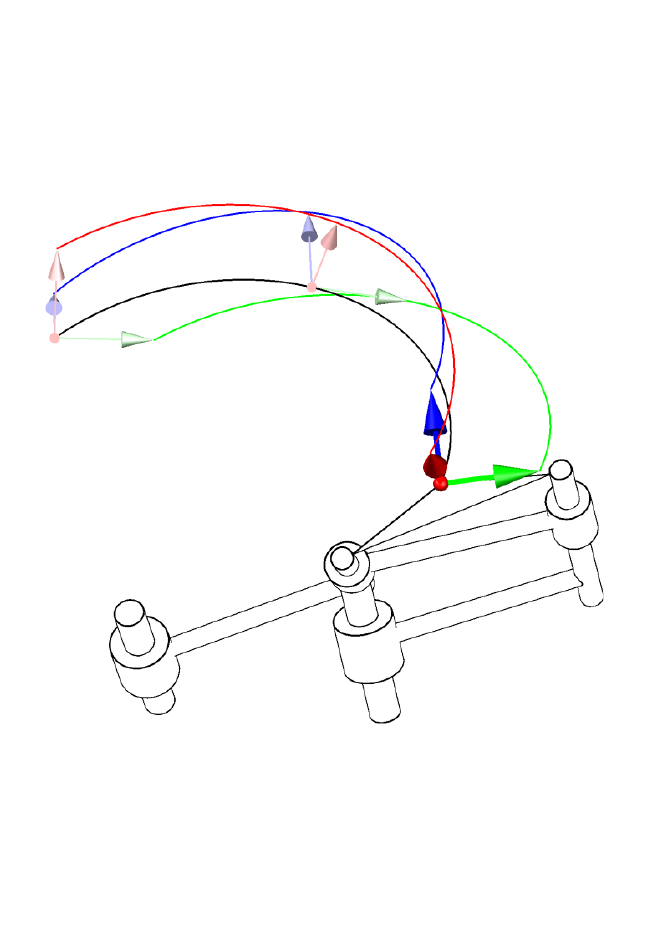}
    \put(5,81){\contour{white}{$p_0$}}
    \put(55,86){\contour{white}{$p_1$}}
    \put(78,51){\contour{white}{$p_2$}}
    \put(7,22){$h_1$}
    \put(7,53){$h_2$}
    \put(49,25){$k_1$}
    \put(50,48){$k_2$}
  \end{overpic}
  \caption{Three-pose synthesis of a four-bar linkage}
  \label{fig:bennett}
\end{figure}

Let us consider a simple example of a spatial four-bar linkage
(Bennett linkage). Its coupler motion is given by a quadratic motion
polynomial $C$. In the kinematic image space $P^7$, it parametrizes a
conic section on the Study quadric $\SQ$. Conversely, a generic conic
section on $\SQ$ gives rise to generic quadratic motion polynomials
(differing only by admissible reparametrizations) and, via motion
factorization, to a Bennett linkage. Thus, we can synthesize a Bennett
linkage to three prescribed poses $p_0$, $p_1$, $p_2$. These poses are
points on the Study quadric $\SQ$ where they span a plane. We compute
a rational quadratic parametrization $C$ of the intersection conic of
this plane and $\SQ$ and factor it as
$C = (t-h_1)(t-h_2) = (t-k_1)(t-k_2)$. The linear motion polynomials
in these factorizations determine the fixed axes ($h_1$, $k_1$) and the
moving axes ($h_2$, $k_2$), as illustrated in Figure~\ref{fig:bennett}.

\begin{figure}
  \small
  \centering
  \begin{overpic}{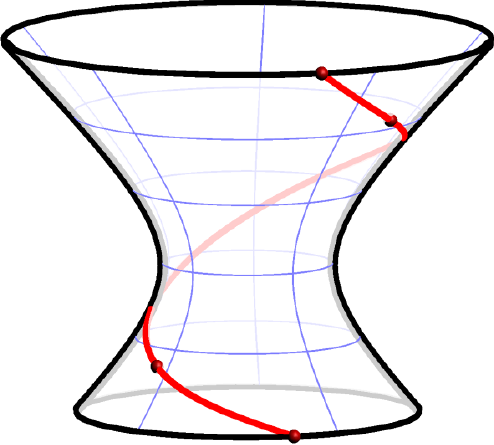}
    \put(54,-5){$p_0$}
    \put(33,20){\contour{white}{$p_1$}}
    \put(77,69){\contour{white}{$p_2$}}
    \put(62,78){\contour{white}{$p_3$}}
  \end{overpic}
  \hfill
  \begin{overpic}{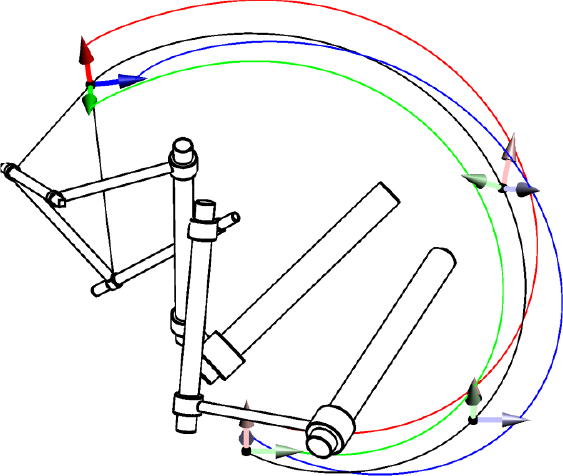}
    \put(60,52){$h_1$}
    \put(30,62){$h_2$}
    \put(6,54){$h_3$}
    \put(78,41){$k_1$}
    \put(38,34){$k_2$}
    \put(10,30){$k_3$}
    \put(38,0){$p_0$}
    \put(87,13){\contour{white}{$p_1$}}
    \put(94,55){\contour{white}{$p_2$}}
    \put(9,70){$p_3$}
  \end{overpic}
  \caption{Four-pose synthesis of a six-bar linkage}
  \label{fig:6R}
\end{figure}

Our synthesis procedure for Bennett linkages is elegant and simple but
it is not the only available method. Other approaches include
\cite{veldkamp67,suh69,tsai73,perez04,brunnthaler05}. It extends,
however, to linkages with more than four joints. We may, for example,
synthesize six-bar linkages to four prescribed poses
\cite{hegedus15:_four_pose_synthesis}. Here, the geometry is slightly
more involved but still understandable via classical results. Four
points $p_0,p_1,p_2,p_3$ in general position on the Study quadric
$\SQ$ span a three-dimensional projective space. If this space
intersects $\SQ$ in a ruled quadric, as in Figure~\ref{fig:6R}, the four
points can be interpolated by two one-parametric families of rational
cubic curves. Each interpolating cubic can be factored and gives rise
to several open chains with three revolute joints that can be combined
to form a closed six-bar linkage. Here are a few further remarks on
this construction:
\begin{itemize}
\item Closed loops of six revolute axes are, in general, rigid. In our
  case, they move because of their special algebraic construction.
\item There exist spatial six-bar linkages that have a one-parametric
  mobility but cannot be constructed by factorization of motion
  polynomials. The relative motions between two links are not all
  rational or, at least, do not have rational components.
\item The above construction is the first viable synthesis procedure for
  spatial six-bar linkages. It may easily be generalized using a Hermite like
  interpolation scheme.
\end{itemize}

\section{Exceptional factorizations and Kempe's Universality Theorem}
\label{sec:exceptional-factorizations}

Theorem~\ref{th:1} shows the existence of finitely many factorizations of
generic motion polynomials. In this section we are concerned with
non-generic situations where the primal part of the motion polynomial
has real factors. In this case, Theorem~\ref{th:1} gives no
information. The following examples demonstrate what can happen:

\begin{example}
  The motion polynomial $C = (t-1)(t - \qj) - \eps((\qi+\qk)t - 2\qk)$
  can be factored as
  \begin{equation*}
    C = (t - 1 - \eps\qi)(t - \qj - \eps\qk)
      = (t - \qj - \eps(\qi+2\qk))(t-1+\eps\qk).
  \end{equation*}
  The polynomial factors $t - 1 - \eps\qi$ and $t - 1 - \eps\qk$
  parametrize, however, translations, not rotations. One may view this
  as a limiting case of the generic situation appearing in
  Theorem~\ref{th:1} in the sense that the two rotations degenerate to
  translations. It turns out that they can still be computed by
  Algorithm~\ref{alg:factor}.
\end{example}

\begin{figure}
  \centering
  \small
  \begin{overpic}{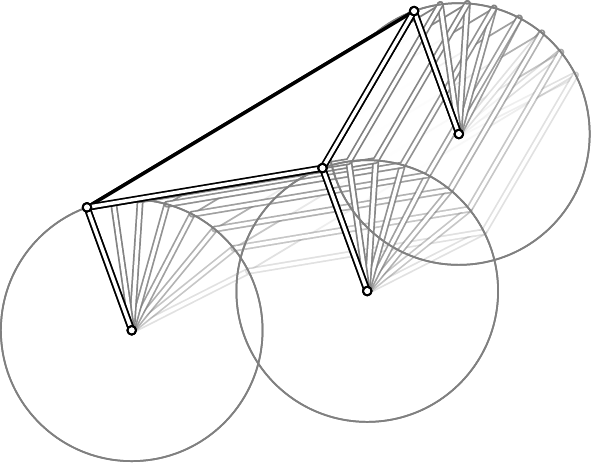}
    \put(14,20){$h''_1$}
    \put(8,46){$h''_2$}
    \put(64,24){$h'_1$}
    \put(49,53){$h'_2$}
    \put(79,53){\contour{white}{$h_1$}}
    \put(63,78){$h_2$}
  \end{overpic}
  \caption{Three different factorizations of a circular translation}
  \label{fig:parallelogram}
\end{figure}
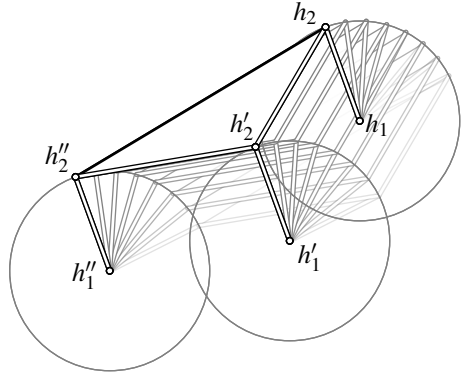

\begin{example}
  \label{ex:3}
  We consider the motion polynomial
  $C = t^2 + 1 + \eps (a\qi + b\qj t)$ with $a, b \in \R$,
  $a, b \ge 0$, $a^2 + b^2 > 0$. It parametrizes the curvilinear
  translation along an ellipse with semi-axis lengths $a$ and $b$. If
  $a \neq b$, a straightforward computation shows that $C$ admits no
  factorization of the form $C = (t - h_1)(t - h_2)$ with linear
  motion polynomials $t-h_1$ and $t-h_2$. If, however, $a = b$
  (circular translation), even infinitely many factorizations exist,
  namely
  \begin{equation}
    \label{eq:5}
    h_1 = \qk - \eps (f \qi + (a + g) \qj),
    \quad
    h_2 = -\qk + \eps (f \qi + g \qj),
  \end{equation}
  with $f,g \in \R$. They have a very clear geometric
  explanation. The motion in question is a circular translation and
  can be generated by infinitely many parallelogram linkages
  (Figure~\ref{fig:parallelogram}). Each leg corresponds to one
  factorization of the form \eqref{eq:5}.
\end{example}

Existence of exceptional situations demonstrate that the factorization
theory over $\DH[t]$ is more complicated but also more interesting
than the theory over $\H[t]$. Moreover, situations with no or
infinitely many factorization arise surprisingly often in engineering
applications of motion polynomial factorization. Many important
rational motions are not amenable to straightforward factorization via
Algorithm~\ref{alg:factor}. Still, it is possible to factor these motions
but at the cost of raising the number of factors (the degree of the
motion polynomial). Recall that for any non-zero polynomial
$R \in \R[t]$ the motion polynomials $C$ and $CR$ parametrize the same
motion. Thus, one may try to find a real polynomial $R$ such that $CR$
admits a factorization.

In engineering applications, revolute joints are preferred over
translational joints. Hence, we focus on factorizations with revolute
joints only. In this case, a necessary requirement is that the motion
parametrized by $C$ is bounded, i.e., all trajectories are bounded
rational curves. We also say that the motion polynomial $C$ itself is
\emph{bounded.} Indeed, we have

\begin{theorem}
  \label{th:2}
  For every bounded motion polynomial $C$ there exists a real
  polynomial $R$ such that $CR$ admits a factorization
  $CR = (t-h_1) \cdots (t-h_m)$ with rotation polynomials $t-h_1$,
  \ldots, $t-h_n$. The degree of $R$ is bounded by the maximal degree
  of a real factor of the primal part of~$C$.
\end{theorem}

Theorem~\ref{th:2} has been proved in \cite{gallet15} for the planar case
and in \cite{li15:_motion_polynomials} for the general case. The
proofs are constructive so that factorizations can be effectively
computed. The main ingredients are variants of the Euclidean algorithm
for polynomial division over the (dual) quaternions and the solution
of quadratic equations with real coefficients over quaternions
\cite{huang02}. It is worth noting that infinitely many factorizations
exist if $\deg R > 0$.

\begin{figure}
  \centering
  \includegraphics[page=33,trim=40 110 70 15,clip]{img/ellipse1}
  \includegraphics{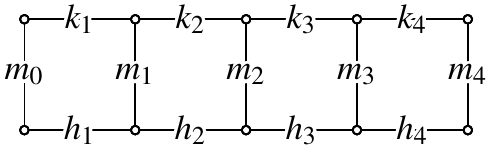}
  \caption{Linkage to generate an elliptic translation and
    corresponding link graph}
  \label{fig:elliptic-translation}
\end{figure}

As an application of Theorem~\ref{th:2}, consider the elliptic translation
appearing in Example~\ref{ex:3}. There exists a quadratic real polynomial
$R$ such that $CR$ admits the factorization
$CR = (t-h_1)(t-h_2)(t-h_3)(t-h_4)$. From this, we construct an
(admittedly complicated) linkage to generate an elliptic translation
(Figure~\ref{fig:elliptic-translation}):
\begin{itemize}
\item The linkage consists of four anti-parallelograms
  ($h_im_ik_im_{i-1}$ for $i\in\{1,2,3,4\}$).
\item Six angles ($\sphericalangle(h_i,m_i,h_{i+1})$,
  $\sphericalangle(k_i,m_i,k_{i+1})$, for $i \in \{1,2,3\}$) are kept
  constant during the motion. In other words, we have a chain of
  anti-parallelograms, where each anti-parallelogram
  ``follows'' its predecessor.
\item The rigid body attached to the connection of $m_4$ and $h_4$
  performs the elliptic translation. The point $h_4$ ``draws'' the
  indicated ellipses.
\end{itemize}
Figure~\ref{fig:elliptic-translation} also shows a more abstract
representation of the same linkage. In this ``link graph'', each
vertex represents a link (a rigid connection between revolute joints)
and each edge represents a joint. Two vertices are connected, if the
corresponding links share a joint. The linkage of Figure~\ref{fig:elliptic-translation} was
constructed by augmenting the given factorization ($h_1$, $h_2$,
$h_3$, $h_4$) with one additional joint $m_0$ and then successively
computing the remaining joints by solving the recursion
\begin{equation}
  \label{eq:6}
  (t-m_{i-1})(t-h_i) = (t-k_i)(t-m_i),
  \quad i \in \{1,2,3,4\}
\end{equation}
with the aid of Algorithm~\ref{alg:factor}. We call this computation of $k_i$
and $m_i$ from $m_{i-1}$ and $h_i$ a \emph{Bennett flip.} This name
comes from the fact that in the spatial case the involved dual
quaternions determine the axes of a Bennett linkage. In the planar
case, a Bennett flip generates anti-parallelograms.

\begin{figure}
  \centering
  \includegraphics{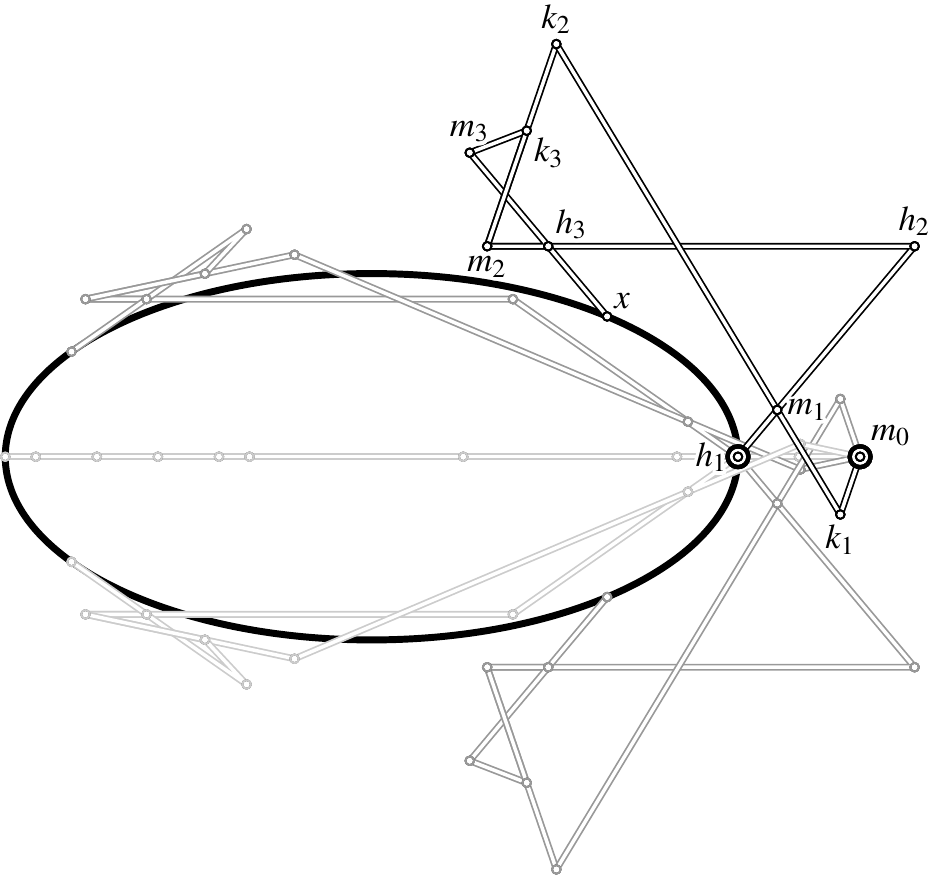}
  \includegraphics{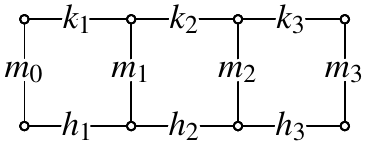}
  \caption{Linkage to draw an ellipse and link graph}
  \label{fig:ellipse2}
\end{figure}

Figure~\ref{fig:ellipse2} presents a refinement of this
construction. Instead of multiplying $C$ with a polynomial
$R \in \R[t]$, we right multiply it with the linear polynomial
$H \coloneqq t-\qk \in \H[t]$ such that the product $CH$ admits a
factorization. Of course, $C$ and $CH$ parametrize different motions
but the fixed points of $H$ (points of the third coordinate axis)
retain their trajectories. Hence, we obtain a linkage that is capable
of drawing the prescribed ellipse. The triples of collinear joints
come from the fact that the constant angles in the linkage of
Figure~\ref{fig:elliptic-translation} are straight for the linkage of
Figure~\ref{fig:ellipse2}. The motion of the link connecting $m_3$ and
$h_3$ is no longer an elliptic translation and generic trajectories
are of degree greater than two.

It is a yet unpublished result that one can always find a polynomial
$H \in \H[t]$ such that $CH$ admits a factorization. Using the Bennett
flip technique, one can then construct a (in general) spatial linkage
with prescribed trajectories. The anti-parallelograms of
Figure~\ref{fig:ellipse2} become Bennett linkages and the additional
freedom we gain by multiplying with a quaternion polynomial
$H \in \H[t]$ instead of a real polynomial $R \in \R[t]$ can be used
to reduce the number of anti-parallelograms.

The constructions of this section are not restricted to ellipses but
can be generalized to arbitrary rational space curves. First let $D$
denote a rational curve of degree $d$ in the Euclidean 3-space. Then
consider a motion polynomial $C$ that parametrizes a rational motion
with trajectory $D$, for example the curvilinear translation along $D$.
Using the algorithm of Theorem~\ref{th:2} we arrive at a factorization of
$CR$ as $(t-h_1)\ldots (t-h_m)$. This gives us an open chain of $m$
revolute joints representing our motion polynomial $C$. Using Bennett
flips as in \eqref{eq:6}, we can construct linkages with only revolute
joints that draw arbitrary (bounded) rational curves.

This is not a new result. By a celebrated theorem of mechanism science
(Kempe's Universality Theorem \cite{kempe76,demaine07}), any bounded
portion of an algebraic curve can be drawn by a linkage. Our approach
shows that in the rational case the number of necessary links and
joints decreases dramatically. The asymptotic bound for rational
curves is linear in the curve degree $n$ while the currently best
known bound for general algebraic space curves is cubic
\cite{abbott08}. Moreover, our general construction behaves quite well
in important low-degree cases. The linkage in Figure~\ref{fig:ellipse2}
has only ten joints while the ellipse linkage in Kempe's construction
requires as much as 235 joints. Often it is possible to further reduce
the number of joints so that engineering applications come into
reach. For example, only seven joints are necessary to generate the
so-called Darboux motion, a spatial motion where all trajectories are
ellipses in non-parallel planes
\cite{li15:_darboux_linkages,li14:_straight}.

\section{Conclusion}
\label{sec:conclusion}

Motion polynomial factorization is an algebraic theory with surprising
relations to mechanisms science. The interplay between both
disciplines is beneficial to both. Algebra can provide solutions to
hitherto inaccessible engineering problems and requirements of
applications led to interesting algebraic questions. Our current work
focuses on both, the details of a version of Kempe's Universality
Theorem for rational space curves with emphasis on a low number of
links and joints and further applications of motion polynomial
factorization to engineering problems.

\section*{Acknowledgments}
\label{sec:acknowledgements}

This work was supported by the Austrian Science Fund (FWF): P~26607
(Algebraic Methods in Kinematics: Motion Factorization and Bond
Theory).

 \bibliographystyle{imn}
\par\noindent
\providecommand{\bysame}{\leavevmode\hbox to3em{\hrulefill}\thinspace}
\providecommand{\MR}{\relax\ifhmode\unskip\space\fi MR }
\providecommand{\MRhref}[2]{  \href{http://www.ams.org/mathscinet-getitem?mr=#1}{#2}
}
\providecommand{\href}[2]{#2}

\end{document}